\documentclass[11pt]{article}
\usepackage{latexsym}
\usepackage[all]{xy}
\usepackage{amsmath}
\usepackage{amssymb}
\usepackage{amsfonts}
\usepackage[applemac]{inputenc}

\newcommand  {\lb} {[\![}
\newcommand  {\rb} {]\!]}

\newtheorem{thm}{Theorem}[section]
\newtheorem{prop}[thm]{Proposition}

\newtheorem{df}[thm]{Definition}
\newtheorem{cor}[thm]{Corollary}

\newtheorem{rmk}[thm]{Remark}

\begin{document}

\title{\textbf{Derived critical loci I  - Basics}} 

\author{Gabriele Vezzosi\\
\small{Dipartimento di Sistemi ed Informatica} \\ 
\small{Universit\`a di Firenze}\\
\small{Italy}}

\date{\textsf{Notes} -- September 2011}

\maketitle

\tableofcontents

\section{Introduction}
We will quickly explore the derived geometry of zero loci of sections of vector bundles, with particular emphasis on derived critical loci. In particular we will single out many of the derived geometric structures carried by derived critical loci: the homotopy \emph{Batalin-Vilkovisky} structure, the \emph{action of the $2$-monoid} of the self-intersection of the zero section, and the \emph{derived symplectic structure} of degree $-1$, and show how this structure exists, more generally, on \emph{derived lagrangian intersections} inside a symplectic manifold. These are just applications of a small part of a much larger project investigating quantization of derived moduli spaces (\cite{ptvv}).\\

These pages were motivated by a series of lectures - on the \emph{formal} derived aspects of the same topics treated here - by Kevin Costello at DennisFest 2011, Stony Brook. The idea was to describe those global structures on a derived critical locus, that after passing to the formal completion (i.e. to the associated \emph{formal moduli problem} in the sense of \cite{lu}) could recover Costello-Gwilliam's set-up in \cite{cg}. \\
The aim of these notes is more at giving a roadmap and an overview of all the constructions than at giving detailed proofs, which will probably appear inside the more comprehensive \cite{ptvv}. In this first part of the notes, we describe the very basics, while part II (\cite{V}) will deal with the explicit induced tame $E_1$-deformations of the periodic derived monoidal dg-category as being developed in \cite{ptvv}.\\

I wish to thank Kevin Costello, Bertrand To\"en and Boris Tsygan for many interesting discussions on these notes' topics.\\

\noindent \textbf{Notations and related warnings -} We work over a base field $k$ of characteristic $0$. The other notations are standard or reminded otherwise. Our reference for derived algebraic geometry will be \cite{hagII}. In particular, $\mathbf{dSt}_{k}$ is the model category of derived stacks over $k$, as defined in \cite[2.2.3]{hagII}. \\ As Kevin Costello pointed out to me, there are different notions of Gerstenhaber and Batalin-Vilkovisky structures in the literature, so my choice of terminology might not be standard. In particular, what we call a dg Gerstenhaber bracket has degree 1, and therefore coincides with what Jacob Lurie calls a $P_0$ bracket. And we use the term dg-Batalin-Vilkovisky (dgBV) algebra according to the physics literature, while Costello calls it a Beilinson-Drinfel'd (BD) algebra. All the structures are explicitly defined in the text, so we hope our terminology will not cause insormountable problems to the reader.

\section{Koszul complexes and derived zero loci of sections of vector bundles}
Since there are tons of complexes named after Koszul in the literature, we give here the definitions we'll use and recall some more or less well known related facts. \\

\subsection{Affine case}
Let $R$ be a commutative $k$-algebra, and $P$ a projective $R$-module of finite type. Let $S:=\mathrm{Sym}_{R}(P^{\vee})$ the symmetric algebra on the $R$-dual $P^{\vee}$: we will \emph{always} disregard its internal grading. $S$ is a commutative $R$-algebra. Let $\bigwedge^{\bullet}P^{\vee}$ be the exterior algebra of $P^\vee$ as an $R$-module. Consider the nonpositively graded $S$-module $S\otimes_{R}\bigwedge^{\bullet}P^{\vee}$ graded by $(S\otimes_{R}\bigwedge^{\bullet}P^{\vee})_m := S\otimes_{R}\bigwedge^{-m}P^{\vee}$ for $m\leq 0$. This is a graded $S$-modules which is degreewise projective over $S$. One has the following further structures on $S\otimes_{R}\bigwedge^{\bullet}P^{\vee}$

\begin{itemize}
\item An obvious $R$\emph{-augmentation} $S\otimes_{R}\bigwedge^{\bullet}P^{\vee} \rightarrow (S\otimes_{R}\bigwedge^{\bullet}P^{\vee})_0 =S \rightarrow R$ coming from the canonical (\emph{zero section}) augmentation of $S$.
\item Since $\bigwedge^{\bullet}P^{\vee}$ is a graded commutative $R$-algebra, $S\otimes_{R}\bigwedge^{\bullet}P^{\vee}$ is a \emph{graded commutative $S$-algebra}.
\item There is a natural degree $1$ \emph{differential} on the graded $S$-module $S\otimes_{R}\bigwedge^{\bullet}P^{\vee}$ induced by contraction and the canonical map $$h:R \longrightarrow \mathrm{Hom}_{R}(P,P) \simeq P^{\vee}\otimes_R P.$$ More precisely, if $h(1)=\sum_{j}\alpha_j \otimes x_j$, then 
$$ d: S\otimes_{R}\bigwedge^{n+1}P^{\vee} \longrightarrow S\otimes_{R}\bigwedge^{n}P^{\vee}$$ acts as $$d(a\otimes (\beta_1 \wedge \ldots \wedge \beta_{n+1}))= \sum_{j} a\cdot \alpha_j \otimes \sum_k (-1)^{k}\beta_{k}(x_j)(\beta_1 \wedge \ldots \wedge \check{\beta_k} \wedge \ldots \wedge \beta_{n+1}).$$ This is obviously $S$-linear.
\item Together with such a $d$, $S\otimes_{R}\bigwedge^{\bullet}P^{\vee}$ becomes a \emph{commutative differential non-positively graded} (cdga) over $S$.
\end{itemize}

\begin{df} The $S$-cdga $\mathrm{K}(R;P):=(S\otimes_{R}\bigwedge^{\bullet}P^{\vee}, d)$ is called the \emph{fancy Koszul cdga} of the pair $(R,P)$.
\end{df}

The following result is well known (and easy to verify, or see e.g. \cite{bgs}, \cite[2.3]{mr})

\begin{prop} The cohomology of the fancy Koszul cdga $\mathrm{K}(R;P)$ is zero in degrees $< 0$, and $H^0(\mathrm{K}(R;P))\simeq R$.
\end{prop}

We will be interested in the following translation of the previous result

\begin{cor}\label{fancyres} The augmentation map $\mathrm{K}(R;P) \longrightarrow R$ is a cofibrant resolution of $R$ in the model category $\mathbf{cdga}_{S}$ of cdga's over $S$.
\end{cor}

We now explain the relation of the fancy Koszul cdga to the usual Koszul complex associated to an element $s\in P$, and interpret this relation geometrically.\\ Let $$\mathrm{K}(R,P;s) :=(\bigwedge^{\bullet} P^{\vee}, d_s)$$ be the usual non-positively graded Koszul cochain complex, whose differential $d_s$ is induced by contraction along $s$ (see \cite{bah}). Together with the exterior product $\mathrm{K}(R,P;s)$ is a cdga over $R$.\\ The choice of an element $s\in P$ induces a map of commutative algebras $\varphi_s: S\rightarrow R$ (corresponding to the evaluation-at-$s$ map of $R$-modules $P^{\vee} \rightarrow R$). We will denote $R$ with the corresponding $S$-algebra structure by $R_s$; $R_s$ is a cdga over $S$ (concentrated in degree $0$). We can therefore form the tensor product of cdga's $R_s \otimes_S \mathrm{K}(R;P)$. This is an $S$-cdga, whose underlying graded $S$-module is $$R_s \otimes_S \mathrm{K}(R;P)= R_s \otimes_S (S\otimes_R \bigwedge^{\bullet}P^{\vee})\simeq \bigwedge^{\bullet}P^{\vee},$$ where $\bigwedge^{\bullet}P^{\vee}$ is viewed as an $S$-module via the composite $k$-algebra morphism $$\xymatrix{S \ar[r]^-{\varphi_s} & R \ar[r]^-{\textrm{can}} & \bigwedge^{\bullet}P^{\vee}}.$$ It is easy to verify that, under this isomorphism, the induced differential on $R_s \otimes_S \mathrm{K}(R;P)$ becomes exactly $d_s$. In other words, there is an isomorphism of $S$-cdga's (and therefore of $R$-cdga's, via the canonical map $R\rightarrow S$)  $$R_s \otimes_S \mathrm{K}(R;P) \simeq \mathrm{K}(R,P;s).$$

\subsection{Affine derived zero loci} Let $R$, $P$, and $s$ be as above. Let $X:=\mathrm{Spec}\, R$, and $\mathbb{V}(P):= \mathrm{Spec}(S)$. The canonical map of $k$-algebras $R\rightarrow S$, indices a map $\pi:\mathbb{V}(P)\rightarrow X$, making $\mathbb{V}(P)$ a vector bundle on $X$ with $R$-module of sections $P$. The zero section $0: X\rightarrow \mathbb{V}(P)$ corresponds to the natural augmentation $S\rightarrow R$, while the section $s$, corresponding to the element $s \in P$, induces the other augmentation we denoted as $\varphi_s:S\rightarrow S$ before. Now we con consider the homotopy pullback $$\xymatrix{Z^{h}(s) \ar[r] \ar[d] & X \ar[d]^-{0} \\ X \ar[r]_-{s} & \mathbb{V}(P)}$$ in the model category $\mathbf{dSt}_k$ of \'etale derived stacks over $k$ (\cite[2.2.3]{hagII}). The derived affine stack $Z^{h}(s)$ is called the \emph{derived zero locus} of the section $s$. By definition of homotopy pullback in $\mathbf{dSt}_k$,
we may choose any cofibrant replacement of the natural augmentation $S\rightarrow R$ in the category of $S$-cdga's, for example (Cor. \ref{fancyres}) $\mathrm{K}(R;P) \rightarrow R$, and get $$Z^{h}(s) \simeq \mathbb{R}\mathrm{Spec}(R_s \otimes_{S}^{\mathbb{L}} R)\simeq \mathbb{R}\mathrm{Spec}(R_s \otimes_{S} \mathrm{K}(R;P)) \simeq \mathbb{R}\mathrm{Spec}(\mathrm{K}(R,P;s))$$ by the computation above. In other words

\begin{prop} The usual Koszul cdga $\mathrm{K}(R,P;s)$ is the algebra of functions on the \emph{derived zero locus} of the section $s$. 
\end{prop}

\begin{rmk} \emph{The corresponding statement in the $C^{\infty}$-category can be found in \cite[Appendix]{cg}.}
\end{rmk}

\subsection{General - i.e. non-affine - case}
Let $X$ be a scheme over $k$, $E$ a vector bundle on $X$, $s\in \Gamma(X,\mathcal{E})$ a section of $E$, and $$\mathrm{K}(X,E;s):=(\bigwedge^{\bullet} \mathcal{E}^{\vee}, d_s)$$ the usual non-positively graded Koszul cochain complex, whose differential $d_s$ is induced by contraction along $s$. Together with the exterior product $\mathrm{K}(X,E;s)$ is a sheaf of $\mathcal{O}_{X}$-cdga's. 
Define the derived zero locus $Z^{h}(s)$ via the homotopy cartesian diagram (in $\mathbf{dSt}_{k}$) $$\xymatrix{Z^{h}(s) \ar[r] \ar[d] & X \ar[d]^-{s} \\ X \ar[r]_-{0} & E.}$$

The following result is an immediate consequence of the corresponding result in the affine case

\begin{prop} $\mathrm{K}(X,E;s)$ is the cdg-Algebra $\mathcal{O}_{Z^{h}(s)}$ of functions on the \emph{derived zero locus} $Z^{h}(s)$ of the section $s$. 
\end{prop}

\section{The self-intersection 2-monoid and its actions}
Let $X=\mathrm{Spec}\, R$ be a smooth $k$-scheme, and $P$ be a finitely generated projective $R$-module. The derived self-intersection $0^2:=Z^{h}(0)$ of the zero section $0$ of $\pi:\mathbb{V}(P)\rightarrow X$ has, \emph{relative to} $X$, one monoid structure (by ``loop'' composition) and one abelian group structure (by addition of ``loops'', coming from the group structure of $\mathbb{V}(P)$ ). These structures are compatible in the sense that the functor corresponding to $0^2$ factors through the category of group-like Segal monoid objects in simplicial abelian groups. Let's denote this structure by $(0^2, \circ, +)$, and observe that this whole structure acts, relatively to $X$, on the derived zero locus of any section of $\pi:\mathbb{V}(P)\rightarrow X$.\\

\emph{More precisely}, the derived self-intersection $0^2$ of the zero section $0$ of $\pi:\mathbb{V}(P)\rightarrow X$ is $$0^2\simeq\mathbb{R}\mathrm{Spec}(\mathrm{K}(R,P;0))\simeq \mathbb{R}\mathrm{Spec}(\bigwedge^{\bullet} P^{\vee}, 0).$$

\begin{rmk} \emph{Note that if $R=k$ and $P$ is free of rank $1$, then $0^2\simeq\mathbb{R}\mathrm{Spec} (k[\varepsilon])$, where $\deg \,\varepsilon=-1$. Also, if $R$ is arbitrary, and $P:=\Omega_{R/k}^1$, then $0^2$ is the derived spectrum of the dg-algebra of polyvectorfields on $X$ - arranged in non-positive degrees and with zero differential. }
\end{rmk}

The derived Cech nerve of the inclusion $0: X\hookrightarrow \mathbb{V}(P)$ is a Segal groupoid object in the model category of derived stacks \emph{over} $X$, the composition law being given - as usual for derived Cech nerves - by $$\xymatrix{(X\times_{\mathbb{V}(P)}^{h} X) \times_{X}^{h} (X\times_{\mathbb{V}(P)}^{h} X) \ar[r]^-{\sim} & X\times_{\mathbb{V}(P)}^{h} X \times_{\mathbb{V}(P)}^{h} X \ar[r]^-{p_{13}} & X\times_{\mathbb{V}(P)}^{h} X}.$$

The linear structure on $\mathbb{V}(P)$ allows one to define another composition law on $0^2$ by simply pointwise adding the ``loops'' : this corresponds to the natural dg-coalgebra structure on $(\bigwedge^{\bullet} P^{\vee}, 0)$. 

Summing up, we get

\begin{prop}\label{nonaffzerolocus} There is a functor $$Z^{h}(0)^{\bullet \, \bullet}:\Delta^{\textrm{op}} \times \Gamma \longrightarrow \mathbf{dSt}_{k}/X$$ such that 
\begin{itemize} 
\item it prolonges the two previous one-level constructions;
\item it defines a group-like Segal monoid object in commutative Segal monoid objects in $\mathbf{dSt}_{k}/X$;
\item $Z^{h}(0)^{\bullet \, \bullet}$ is well defined in the homotopy category of group-like Segal monoid objects in commutative Segal monoid objects in $\mathbf{dSt}_{k}/X$;
\end{itemize}
\end{prop}

Here, as customary, $\Gamma$ is the category of finite pointed sets and pointed maps. Observe that the $\Gamma$ argument corresponds to addition of loops while the $\Delta^{\textrm{op}}$ one to the composition of loops.\\

\begin{rmk}\emph{The results in Prop. \ref{nonaffzerolocus} are contained in \cite{pre}, in the case where $P$ is the free $k$-module of rank $1$ (i.e. the trivial line bundle $\mathbb{A}_{X}^{1}$ on $X=\mathrm{Spec} \, k$), and interestingly related to matrix factorizations dg-derived categories.}
\end{rmk}

We will call a group-like Segal monoid object in commutative Segal monoid objects (in $\mathbf{dSt}_{k}/X$), a $\textsf{(G,Comm)}$ $2$\emph{-monoid} (in $\mathbf{dSt}_{k}/X$).\\

\begin{prop} For any section $s$ of $\mathbb{V}(P)\rightarrow X$, there is an action of the $\textsf{(G,Comm)}$ $2$\emph{-monoid} (in $\mathbf{dSt}_{k}/X$) $Z^{h}(0)^{\bullet \, \bullet}$ on $Z^{h}(s)$ over $X$.
\end{prop}

\begin{rmk}\emph{Similar results obviously hold in the non-affine case. We leave to the reader the parallel statements of such results.}
\end{rmk}

\section{Derived critical loci}
Let $M$ be a smooth algebraic variety over $k$, $f\in \mathcal{O}(M)$, and $df$ the induced section of the cotangent bundle of $M$.

\begin{df} The \emph{derived critical locus} $\mathrm{Crit}^{h}(f)$ of $f$ is the derived zero locus $Z^{h}(df)$ of $df$.
\end{df}

As a corollary of Proposition \ref{nonaffzerolocus}, we have

\begin{prop} $\mathrm{K}(M,T^{\vee}M;df)$ is the cdg-Algebra $\mathcal{O}_{\mathrm{Crit}^{h}(f)}$ of functions on the \emph{derived critical locus} $\mathrm{Crit}^{h}(f)$ of the function $f$. 
\end{prop}

In the next Section we will describe some of the derived geometrical structures carried by derived critical loci.

\section{Geometric structures on derived critical loci}

\subsection{$G_{\infty}$-structure}

For Let $X=\mathrm{Spec}\,R$, we will denote by $\mathrm{T}_X$ the $R$-dual module to $\Omega^{1}(X)\equiv \Omega_{R/k}^{1}$. For definitions and basic properties of differential graded (homotopy) Gerstenhaber algebras we refer the reader to \cite{ma,tt,ge}. \\

\begin{prop} Let $X=\mathrm{Spec}\,R$ be a \emph{smooth} affine scheme over $k$, and $\alpha \in \Omega^1(X)$ a \emph{closed} $1$-form on $X$. Then, the usual Koszul cdga $(\mathrm{K}(R, \Omega^1(X); \alpha), d_{\alpha})=(\bigwedge^{\bullet}\mathrm{T}_X, d_{\alpha})$, together with the Schouten bracket $\lb -, - \rb$ is a differential graded commutative Gerstenhaber algebra.
\end{prop}

\noindent \textsf{Proof.} The result follows from the easily verified formula $$\alpha([X,Y])= d_{\alpha} \lb X, Y \rb = \lb d_{\alpha}X, Y \rb + \lb X, d_{\alpha}Y \rb = -Y(\alpha(X))+ X(\alpha(Y))$$ valid for any $X, Y \in \mathrm{T}_X$, $[-,-]$ being the usual Lie bracket of vector fields on $X$. \hfill $\Box$ \\

Since $$Z^{h}(\alpha) \simeq \mathbb{R}\mathrm{Spec} (\bigwedge^{\bullet}\mathrm{T}_X, d_{\alpha})$$ (isomorphism in the homotopy category of $\mathbf{dSt}_k$) - i.e. $(\bigwedge^{\bullet}\mathrm{T}_X, d_{\alpha})$ is a model for the cdga of functions on $Z^{h}(\alpha)$ - by the general $G-G_{\infty}$ principle ($\equiv$ ``strict $G$ on a dg-model $\Rightarrow$ $G_{\infty}$ on any dg-model''), we get

\begin{cor} Let $X=\mathrm{Spec}\,R$ be a \emph{smooth} affine scheme over $k$, $\alpha \in \Omega^1(X)$ a \emph{closed} $1$-form on $X$, and $Z^{h}(\alpha)$ the derived zero locus of $\alpha$. Then the cdga $\mathcal{O}_{Z^{h}(\alpha)}$ of functions on $Z^{h}(\alpha)$ is a dg-$G_{\infty}$-algebra.
\end{cor}

\subsection{$BV_{\infty}$-structure}
If $X$ has a volume form $\textsf{vol}$ (e.g. $X$ is symplectic or Calabi-Yau), then we may consider the differential $\Delta$ induced by the de Rham differential via the contraction-with-$\textsf{vol}$ isomorphism $$i_{\textsf{vol}}: \bigwedge ^{\bullet} T_{X} \longrightarrow \bigwedge ^{n-\bullet} \Omega^{1}_{X}$$ where $n:=\dim X$. Together with the bracket defined above, this gives an explicit $BV$ structure on the dg-model $(\bigwedge^{\bullet}\mathrm{T}_X, d_{\alpha})$ hence a dg-$BV_{\infty}$-structure on $\mathcal{O}_{Z^h(\alpha)}$. 

\subsection{Derived symplectic structures}
We will first define degree $n$ derived symplectic structures on a derived stack, and then show how derived intersections of smooth lagrangians (in particular derived critical loci) carry canonical degree $(-1)$ derived symplectic structures.
 
\subsubsection{Generalities on derived symplectic structures}
Derived symplectic structures are introduced and studied in detail in \cite{ptvv} where, in particular, they are shown to induce $E_{n}$-deformation of the derived dg-category of the underlying derived stack. We will only briefly sketch their definitions here, and show how they naturally arise on derived intersections of lagrangian submanifolds.\\

\begin{df} Let $Y$ be a derived algebraic stack over $k$, $\mathcal{L}Y:= \mathbb{R}\mathrm{HOM}_{\mathbf{dSt}_{k}}(S^1:=B\mathbb{Z}, Y)$ its \emph{derived loop stack}, and $\widehat{\mathcal{L}Y}$ its \emph{derived formal loop stack} (i.e. the derived completion of $\mathcal{L}Y$ along the constant loops $Y$) together with the canonical action of $\mathbb{G}_{m}\ltimes B\mathbb{G}_{a}$.
\begin{itemize}
\item The \emph{Hochschild homology} complex of $Y$ is $\mathrm{HH}(Y):= \mathcal{O}_{\widehat{\mathcal{L}Y}}$, together with the additional internal grading by the $\mathbb{G}_{m}$-action; the corresponding internal degree will be called the \emph{weight}.
\item  The \emph{negative cyclic homology} complex of $Y$ is $\mathrm{HC}^{-}(Y):= \mathcal{O}_{\widehat{\mathcal{L}Y}}^{B\mathbb{G}_{a}}$, together with the additional internal grading by the residual $\mathbb{G}_{m}$-action; the corresponding internal degree will be called the \emph{weight}.
\item Let $i: \mathrm{HC}^{-}(Y) \longrightarrow \mathrm{HH}(Y)$ be the obvious canonical (weight-preserving) map.
\end{itemize}
\end{df}

Note that if $\mathbb{L}_{Y}$ is the cotangent complex of $Y$, an element $\omega \in \mathrm{HH}_{m}(Y)$ corresponds to a map $\omega: \mathcal{O}_{Y} \rightarrow \mathrm{Sym}^{\bullet}(\mathbb{L}_{Y}[1])[-m]$ in $\mathrm{D}(Y)$, and its weight is just its degree in $\mathrm{Sym}^{\bullet}$. \\
Recall from \cite{ill} that the derived functors $\wedge_{\mathcal{O}_{Y}}^{m}$ and $\mathrm{Sym}_{\mathcal{O}_{Y}}^{m}$, $m\geq 0$ are related by $$\wedge_{\mathcal{O}_{Y}}^{m}(P)\simeq \mathrm{Sym}_{\mathcal{O}_{Y}}^{m}(P[1])[-m]$$ for any $\mathcal{O}_{Y}$-Module $P$. In particular, we have 
$$\mathrm{Sym}_{\mathcal{O}_{Y}}^{2}(\mathbb{L}_{Y}[1])\simeq \mathbb{L}_{Y}\wedge \mathbb{L}_{Y} [2].$$
Therefore there is a canonical projection $$\mathrm{Sym}^{\bullet}(\mathbb{L}_{Y}[1]) \longrightarrow \mathbb{L}_{Y}\wedge \mathbb{L}_{Y}[2] $$ to the weight $2$ part. \\
 
We are now ready to define generalized or derived symplectic structures.

\begin{df} Let $Y$ be a derived algebraic stack over $k$, $\mathbb{L}_{Y}$ its cotangent complex, and $\mathbb{T}_{Y}=:\mathbb{L}^{\vee}_{Y}$ its tangent complex. A \emph{derived symplectic structure of degree} $n \in \mathbb{Z}$ on $Y$ is an element $\omega \in \mathrm{HC}_{2-n}^{-}(Y)$ of weight $2$ such that its weight $2$ projection $$\overline{\omega}:\mathcal{O}_{Y} \longrightarrow \mathrm{Sym}^{\bullet}(\mathbb{L}_{Y}[1])[n-2] \longrightarrow \mathbb{L}_{Y}\wedge \mathbb{L}_{Y} [n]$$ induces, by adjunction, an isomorphism $\mathbb{T}_{Y}\simeq \mathbb{L}_{Y}[n]$ (non-degeneracy condition).
\end{df}

Intuitively, a derived symplectic form is therefore a non-degenerate map $\mathbb{T}_{Y}\wedge \mathbb{T}_{Y} \rightarrow \mathcal{O}_{Y}[n]$ which is $B\mathbb{G}_{a}$-equivariant (i.e. it lifts to $\mathrm{HC}_{2-n}^{-}(Y)$, and this is a datum).\\

\begin{rmk} \emph{If $Y$ is a (quasi-smooth) derived Deligne-Mumford stack and $\omega \in \mathrm{HC}_{3}^{-}(Y)$ is a $(-1)$ derived symplectic structure on $Y$, let us consider the following construction. By definition of derived exterior and symmetric powers of complexes,
we have, for any complex $P$ on $Y$, a canonical isomorphism $$P \wedge P \simeq \mathrm{Sym}^{2}(P[1])[-2],$$ and taking $P= \mathbb{L}_{Y}[-1]$ ($\mathbb{L}_{Y}$ being the cotangent complex of the derived stack $Y$), we get a canonical isomorphism $$\textrm{can}: \mathbb{L}_{Y}[-1] \wedge \mathbb{L}_{Y}[-1] \simeq \mathrm{Sym}^{2}(\mathbb{L}_{Y})[-2].$$
Since our derived stack $Y$ is endowed with a $(-1)$ derived symplectic form, its image $\overline{\omega} \in \mathrm{HH}_{3}(Y)$ yields a morphism  
$$\hat{\omega} : \mathbb{T}_{Y} \wedge \mathbb{T}_{Y} \longrightarrow \mathcal{O}_{Y}[-1]$$ such that the map $$\varphi_{\omega}: \mathbb{T}_{Y} \longrightarrow \mathbb{L}_{Y}[-1],$$ induced by adjunction, is an isomorphism. Therefore we may consider the ($-2$-shifted) associated Poisson bivector, i.e. the  composition $$\xymatrix{s_{\omega}: \mathrm{Sym}^{2}(\mathbb{L}_{Y})[-2] \ar[r]^-{\textrm{can}^{-1}} & \mathbb{L}_{Y}[-1] \wedge  \mathbb{L}_{Y}[-1] \ar[rr]^-{f_{\omega} \wedge f_{\omega}} & & \mathbb{T}_{Y} \wedge \mathbb{T}_{Y}  \ar[r]^-{\overline{\omega}} & \mathcal{O}_{Y}[-1].}$$ If $j:\mathrm{t}_{0}(Y) \hookrightarrow Y$ denotes the closed immersion of the truncation, then the canonical map $$j^* \mathbb{L}_{Y} \longrightarrow \mathbb{L}_{\mathrm{t}_{0}(Y)} $$ together with the $2$-shifted restricted map $$j^{*}(s_{\omega}[2]): \mathrm{Sym}^{2}(j^* \mathbb{L}_{Y}) \longrightarrow \mathcal{O}_{\mathrm{t}_{0}(Y)}[1]$$ define a \emph{symmetric perfect obstruction theory} - in the sense of \cite{befan} - on $\mathrm{t}_{0}(Y)$.  \\
My present guess is that \emph{all} known symmetric obstruction theories arise this way, i.e. are induced as explained above from a $(-1)$ derived symplectic structure on some (quasi-smooth) derived extension $Y$ of the given DM stack on which they are defined. In particular this would imply that the known ones have indeed more 'structure', i.e. their lift to derived stack corresponds to a \emph{closed}  element in $\mathrm{HH}_{3}(Y)$ - i.e. an  element admitting a lift to $\mathrm{HC}_{3}^{-}(Y)$. And it might be the case that a symmetric obstruction theory that is induced by a $(-1)$ derived symplectic form, is \'etale locally isomorphic to the canonical one existing on 
the derived zero locus of a \emph{closed} $1$-form on a smooth scheme - instead of just an \emph{almost-closed} $1$-form, as in the case of a general symmetric obstruction theory, see \cite{befan}. }
\end{rmk}

Heuristically, a derived symplectic structure of degree $n\geq0$ (resp. $n\leq 0$) on $Y$ induces a \emph{quantization} of $Y$, i.e. a deformation of the derived dg-category $\mathrm{D}(Y)$ (resp. the periodic derived dg-category $\mathrm{D}^{\mathbb{Z}/2}(Y)$) as an $E_n$-monoidal dg-category (resp. as an $E_{-n}$-monoidal dg-category). This is the result that makes the notion of derived symplectic form geometrically interesting. We address the reader to \cite{ptvv} for detailed statements and proofs, while we will come back on this more in notes' style in \cite{V}. \\

\begin{rmk}\emph{Note that $A:=(\bigwedge^{\bullet}\mathrm{T}_X, d_{\alpha})$ is not cofibrant as a cdga over $k$, so its $\Omega_{A/k}^2$ is not directly related to $\mathbb{L}_{Z^h(\alpha)}\wedge \mathbb{L}_{Z^h(\alpha)}\equiv \mathbb{L}_{A/k}\wedge \mathbb{L}_{A/k}$, hence we cannot exhibit an explicit map $A[1]\rightarrow \Omega_{A/k}^2$ in $D(A)$ realizing $\omega$. It is however possible, in the affine case, to compute explicitly the element in $HC_{3}^-$ corresponding to $\omega$, at least in the case $X$ is affine - see \cite{ptvv}. }
\end{rmk}

\subsubsection{Symplectic structure on derived intersections of lagrangian subvarieties}\label{lagra}
Let $(Y, \omega_{0})$ a smooth algebraic symplectic variety, $L_1$ and $L_2$ smooth lagrangian subvarieties in $X$. We will show that the derived intersection $Z:=L_1 \times_{X}^h L_2$ has a canonical $(-1)$-symplectic structure, and therefore by \cite{ptvv} the derived dg-category $D^{\mathbb{Z}/2}(Z)$ has a canonical deformation $D_{\omega_{0}}^{\mathbb{Z}/2}(Z)$ as a monoidal dg-category. \\

Consider the homotopy cartesian diagram $$\xymatrix{Z \ar[r]^-{j_1} \ar[d]_-{j_2} & L_1 \ar[d]^-{i_1} \\ L_{2} \ar[r]_-{i_2} & Y}$$ By definition of homotopy cartesian diagram we have a suitably functorial homotopy $H: j_1^{*}i_1^{*} \sim j_2^{*}i_2^{*}$. By applying this to $\omega_{0}$, and using that each $L_{i}$ is lagrangian, we get a canonical self-homotopy $h_{\omega_{0}}$ of the \emph{zero} map $$j_1^{*}i_1^{*}(T_{Y}\wedge T_{Y}) \longrightarrow \mathcal{O}_{Z}.$$ By composing with the canonical derivative map $\mathbb{T}_{Z}\wedge \mathbb{T}_{Z} \rightarrow j_1^{*}i_1^{*}(T_{Y}\wedge T_{Y})$, we get then a self-homotopy of the zero map $\mathbb{T}_{Z}\wedge \mathbb{T}_{Z} \longrightarrow \mathcal{O}_{Z}$, i.e. a map $$\omega:\mathbb{T}_{Z}\wedge \mathbb{T}_{Z} \longrightarrow \mathcal{O}_{Z}[-1].$$ Hence $\omega$ is an element of $\mathrm{HH}_{3}(Z)$. We let the reader verify that this map lifts canonically to $\mathrm{HC}_{3}^{-}(Z)$ (i.e. it is $B\mathbb{G}_{a}$-equivariant, in other words 'closed') and that the non-degeneracy of $\omega_{0}$ implies that of $\omega$. In other words, $\omega$ is in fact a $(-1)$-symplectic form on the derived intersection $Z:=L_1 \times_{X}^h L_2$.\\

If $Y=T^{*}M$, with $M$ a smooth algebraic variety, then any closed $1$-form $\alpha$ on $M$ defines a Lagrangian subvariety $L_{\alpha}$ in $T^*M$, therefore the derived critical locus $$\mathrm{Crit}^{h}(f)\equiv Z^{h}(df):=L_{df}\times_{T^*M}^h L_{0},$$ admits a derived $(-1)$-symplectic form. \\

\begin{rmk}\emph{The $BV_{\infty}$ structure, the $(-1)$-symplectic structure, and the action of the $\textsf{(G,Comm)}$ $2$-monoid $Z^{h}(0)^{\bullet \, \bullet}$ on a derived critical locus are suitably \emph{compatible}: we address the reader to \cite{ptvv} for details.}
\end{rmk}


\begin{thebibliography}{90}

\bibitem[Be-Fa]{befan} K. Behrend, B. Fantechi, \emph{Symmetric obstruction theories and Hilbert schemes of points on threefolds}, Algebra \& Number Theory, 2 (2008), 313-345.
\bibitem[Be-Gi-So]{bgs} A. Beilinson, V. Ginzburg, W. Soergel, \emph{Koszul duality patterns in representation theory}, JAMS 9 (1996), 473–527.
\bibitem[Bou-AH]{bah} N. Bourbaki, \emph{Algèbre - Ch. 10 - Algèbre homologique}, Springer Verlag, 2007.
\bibitem[Co-Gw]{cg} K. Costello, O. Gwilliam, \emph{Factorization algebras in perturbative quantum field theory}, Preprint draft, 2011.
\bibitem[Ge]{ge} E. Getzler, \emph{Batalin-Vilkovisky algebras and two-dimensional topological field theories}, Comm. Math. Phys. 159 (1994), 265-285. 
\bibitem[Ill]{ill} L. Illusie, \emph{Complexe cotangent et d\'eformations I}, Lecture Notes in Mathematics \textbf{239}, Springer Verlag, Berlin, 1971.
\bibitem[Lu]{lu} J. Lurie, \emph{Formal Moduli Problems}, available at http://www.math.harvard.edu/$\sim$lurie as DAG-X.
\bibitem[Ma]{ma} M. Manetti, \emph{Lectures on deformations of complex manifolds}, Rendiconti di Matematica 24, (2004), 1-183 (and arxiv:math/0507286).
\bibitem[Me]{me} Merkulov, \emph{Wheeled pro(p)file of Batalin-Vilkovisky formalism}, Commun. Math. Phys. 295, No.3 (2010) , 585-638.
\bibitem[Mi-Ri]{mr} I. Mirkovic, S. Riche, \emph{Linear Koszul duality}, Comp. Math. 146 (2010), 233-258.
\bibitem[Ta-Ts]{tt} D. Tamarkin, B. Tsygan, \emph{Noncommutative differential calculus, homotopy BV algebras and formality conjectures}, Methods Funct. Anal. Topology 6 (2000), no. 2, 85-100.
\bibitem[Pa-To-Va-Ve]{ptvv} T. Pantev, B. To\"en, M. Vaqui\'e, G. Vezzosi, \textit{Quantization and derived moduli spaces}, Preprint, in preparation.
\bibitem[Pr]{pre} A. Preygel, \emph{Thom-Sebastiani \& Duality for Matrix Factorizations}, Preprint arXiv:1101.5834.
\bibitem[HAG-II]{hagII} B. To\"en, G. Vezzosi, \textit{Homotopical algebraic geometry II: Geometric stacks
and applications}, Mem. Amer. Math. Soc.  193  (2008),  no. 902, x+224 pp.
\bibitem[To-Ve]{tv} B. To\"en, G. Vezzosi, \textit{Alg\`ebres simpliciales $S^{1}$-\'equivariantes, th\'eorie de de Rham et th\'eor\`emes HKR multiplicatifs}, Compositio Math., in press.
\bibitem[Ve]{V} G. Vezzosi, \textit{Derived critical loci II - quantization}, Notes, in preparation.





\end{thebibliography}
\end{document}